# Generalization of Euler-Lagrange Equations to Find Min-max Optimal Solution of Uncertain Systems

Farid Sheikholeslam, R. Doosthoseyni

Department of Electrical and Computer Engineering, Isfahan University of Technology, Isfahan, IRAN

**Abstract**

In this paper, calculus of variation methods are generalized to find min-max optimal solution of uncertain dynamical systems with uncertain or certain cost. First, a new form of Euler-Lagrange conditions for uncertain systems is presented. Then several cases are indicated where final condition can be specified or free. Also necessary conditions are introduced to existence of min-max optimal solution of the uncertain systems. Finally, efficiency of the proposed method is verified through some examples.

## 1. Introduction

THERE are some papers that use Euler-Lagrange equations in a new way. A variational analysis of neutral functional-differential inclusions and obtain new necessary optimality conditions of both Euler--Lagrange and Hamiltonian types is used in [1]. Bounded slope condition on the boundary values of a minimization problem for a functional of the gradient and its application of this result to prove the validity of the Euler Lagrange equation presented in [2]. An Euler-Lagrange inclusion for optimal problems is used in [3]. Direct method to obtain necessary optimality conditions in a refined Euler–Lagrange form without standard convexity assumptions is presented in [4]. Formulation of Euler–Lagrange inclusion for fractional problems

---

*E-mail addresses :* sheikh@cc.iut.ac.ir (Farid Sheikholeslam)

is discussed in [5]. In addition many papers exist such that they solved certain problems using calculus of variational approaches such as Jonckheere [6,7]. optimal control of a class of distributed-parameter systems governed by first-order, quasilinear hyperbolic partial differential equations that arise in optimal control problems of many physical systems such as fluids dynamics and elastodynamics is described in [8].

In practice, all of physical systems have uncertainty and engineer need to model systems with uncertainty [9-11]. However, there exists a few papers described optimal control in present of uncertainty. [12] assumed that uncertainty is in the initial state and state equation and maximizes the cost which controller is attempting to minimize. The min-max approach to uncertainty has also been invested in [13-16]. A Riccati equation approach to the optimal cost control of uncertain systems with structured uncertainty consider in [17]. Also [18] present an adaptive extremum seeking control of nonlinear dynamic systems with parametric uncertainties. They proposed adaptive extremum seeking controller by using inverse optimal to minimize a meaningful cost function. Optimal control of the force developed in an automotive restraint system during a frontal impact study in [19].

In some papers min-max approach is used in min-max predictive control [20]. Min-max approach in ultimate periodicity of orbits has been invested in [21]. They use pure min-max function and conditional redundancy to provide a simple alternate proof to the ultimate periodicity theorem. [22] use state feedback stabilization and majorizing achievement of min-max-plus systems and corresponding algorithm. They give some solutions that construct the state feedback function for a given desired eigenvalue. [23] studied a minimax-based search algorithm for finding a global optimal solution. They presented an algorithm for a class of optimal finite-precision controller realization problem with saddle point. However all of the above papers cannot describe optimal min-max solution of nonlinear dynamics in general. However, there is not general solution for min –max solution of nonlinear uncertain dynamics in all of above papers.

In this paper, min-max solution of uncertain system with certain or uncertain cost is described. The approach is to find optimal solution of uncertain systems with functional constraint. This approach generalized Euler Lagrange equation for nonlinear uncertain systems and parametric dynamical equations with or without functional constraint in general. The problem is presented in different cases such that final condition can be specified or free. Also, necessary conditions to existence of min-max solution of the functional are introduced. Efficiency of the proposed method is verified through some examples. In the next section problem formulation is considered. In section 3, min-max solution of uncertain functional is presented. In section 4, necessary conditions for optimal control problem with uncertainty in state equation and cost function are applied. Also in section 5, necessary conditions for optimal control problem with uncertain parameter in initial state are applied. Conclusions are drawn in section 6.

## 2. Problem formulation

Let $x(t)$ be a scalar function with continuous first derivatives. It is desired to find the trajectory $x^*(t)$ for which the uncertain functional

$$J(x,a) = \int_{t_0}^{t_f} g(x(t), \dot{x}(t), a, t) dt \qquad (1)$$

has a relative min-max solution, where $g : R^n \times R^n \times A \times [t_0, t_f] \to R$ is given function. It is assumed that $g$ has continuous first and second partial derivatives with respect to all of its argument. Initial time and state $t_0$ and $x_0$ are fixed and final time and state $t_f$ and $x_f$ are specified or free respect to problem statement. $a \in A \subset R^m$ is an uncertain parameter such that $A$ is known compact set. If parameter $a$ is specified then the problem reduces to the usual calculus of variations problem where it is desired to find conditions satisfying $x(t)$ is a minimizer [24]. However, the min–max solution of the functional can be formulated in such a way that the operation of the maximization is taken over the set of uncertainty and the operation of the minimization is taken over trajectory. Therefore, the problem is to find an admissible optimal min-max solution $(x^*(t), a^*)$ which satisfying [12]

$$\max_{a \in A} J(x^*(t), a) \leq J(x^*(t), a^*) \leq \min_{a \in A} J(x(t), a^*) \text{ for all } x(t) \in R^n \text{ and } a \in A. \qquad (2)$$

It may note that if the optimal min-max solution $(x^*(t), a^*)$ has a variation in direction of $x(t)$ then the functional increases. On the other hand, if the optimal solution has a variation in direction of $a$ then the functional decreases. It means that, the optimal solution of functional is a saddle point. So, the min-max solution $(x^*(t), a^*)$ is an extremum (saddle point) of the functional.

Here, it is assumed that a saddle point exists, i.e., there is an uncertain parameter $a^*$ maximizing functional, which the trajectory $x^*(t)$ tries to minimize. In the following section necessary conditions to find min-max solution of the functional are presented when final state and final time is specified or free.

## 3. Min-max solution of uncertain functional

In this section, min-max solution of the functional (1) with specified or free boundary conditions are discussed in four cases. All of the theorems in this section follow standard method in calculus of variation and proves of them are similar to the same cases in [24].

### A. Specified final time and final point

The following theorem, which is proved in Appendix I, introduces necessary conditions for the problem with fixed final time and state.

*Theorem 1*: consider uncertain functional (1) subject to

$$x(t_0) = x_0, x(t_f) = x_f, t \in [t_0, t_f], a \in A. \tag{3}$$

Suppose $J(x(t), a)$ has a min-max solution with assumed conditions and let $x^*(t)$ be admissible. If there exists a $a^* \in A$ such that

1) $\dfrac{\partial g}{\partial x}(x^*(t), \dot{x}^*(t), a^*, t) - \dfrac{d}{dt}[\dfrac{\partial g}{\partial \dot{x}}(x^*(t), \dot{x}^*(t), a^*, t)] = 0$

2) $\displaystyle\int_{t_0}^{t_f} \dfrac{\partial g}{\partial a}(x^*(t), \dot{x}^*(t), a^*, t) dt = 0$

3) $\displaystyle\int_{t_0}^{t_f} \dfrac{\partial^2 g}{\partial a^2}(x^*(t), \dot{x}^*(t), a^*, t) dt < 0$

4) $\dfrac{\partial^2 g}{\partial x^2}(x^*(t), \dot{x}^*(t), a^*, t) > 0$

5) $\dfrac{\partial^2 g}{\partial \dot{x}^2}(x^*(t), \dot{x}^*(t), a^*, t) > 0$

6) $x^*(t_0) = x_0, x^*(t_f) = x_f$

then $x^*(t)$ is a min-max solution of (1).

Theorem 1 suggests min-max optimal solution trajectory when there exists a min-max solution. The following example illustrates min-max solution scheme and its theoretic aspect.

*Example 1*: suppose that

$$J(x(t), a) = \int_0^1 \left(x^2(t) + \dot{x}^2(t) - a^2 - 2ax(t)\right) dt, \quad x(0) = 1, x(1) = 1.$$

Note that the integrand doesn't contain $t$ explicitly and uncertain parameter $a$ belongs to $\mathbb{R}$. Thus, from condition 1 of theorem 1 we have

$$\ddot{x} - x + a = 0$$

then

$$x^*(t) = Ae^t + Be^{-t} + a$$

To satisfy the boundary conditions $x(0) = 1, x(1) = 1$, parameters $A$ and $B$ obtain as

$$A = \dfrac{(1-a)}{(e+1)}, B = \dfrac{(1-a)e}{(e+1)}.$$

Thus,

$$x^*(t) = \dfrac{(1-a)}{(e+1)} e^t + \dfrac{(1-a)e}{(e+1)} e^{-t} + a.$$

On the other hand,

$$\dfrac{\partial g}{\partial a}(x(t), \dot{x}(t), a, t) = -2(x+a).$$

Therefore, condition 2 satisfies if

$$0 = \int_0^1 \frac{\partial g}{\partial a}(x^*(t), \dot{x}^*(t), a^*, t)dt = \frac{-4a - 2e + 2}{e + 1}$$

and hence, it follows that

$$a^* = \frac{-e + 1}{2}.$$

Thus, the extremum solution is

$$x^*(t) = \frac{1}{2}e^t + \frac{e}{2}e^{-t} + \frac{-e + 1}{2}$$

and

$$J(x^*, a^*) = e - 1.$$

Also,

$$\int_0^1 \frac{\partial^2 g}{\partial a^2}(x^*(t), \dot{x}^*(t), a^*, t)dt = -2 < 0,$$

$$\frac{\partial^2 g}{\partial \dot{x}^2}(x^*(t), \dot{x}^*(t), a^*, t) = 2 > 0$$

and

$$\frac{\partial^2 g}{\partial x^2}(x^*(t), \dot{x}^*(t), a^*, t) = 2 > 0.$$

So, all conditions are satisfied and the above extremum solution is a min-max solution of the functional.

*B. Specified final time and free final point*

In this case, initial boundary condition and final time are specified but final boundary condition is free. Necessary conditions to find min-max solution of (1) are presented in theorem 2 (proved in Appendix II).

*Theorem 2*: consider uncertain functional (1) subject to

$$x(t_0) = x_0, x(t_f) = \text{free}, \ t_f \text{ is specified}, \ t \in [t_0, t_f], a \in A. \tag{4}$$

Suppose $J(x(t), a)$ has a min-max solution with assumed conditions and let $x^*(t)$ be admissible. If there exists a $a^* \in A$ such that

1) $\frac{\partial g}{\partial x}(x^*(t), \dot{x}^*(t), a^*, t) - \frac{d}{dt}[\frac{\partial g}{\partial \dot{x}}(x^*(t), \dot{x}^*(t), a^*, t)] = 0$

2) $\int_{t_0}^{t_f} \frac{\partial g}{\partial a}(x^*(t), \dot{x}^*(t), a^*, t) dt = 0$

3) $\int_{t_0}^{t_f} \frac{\partial^2 g}{\partial a^2}(x^*(t), \dot{x}^*(t), a^*, t) dt < 0$

4) $\frac{\partial^2 g}{\partial x^2}(x^*(t), \dot{x}^*(t), a^*, t) > 0$

5) $\frac{\partial^2 g}{\partial \dot{x}^2}(x^*(t), \dot{x}^*(t), a^*, t) > 0$

6) $\frac{\partial g}{\partial \dot{x}}(x^*(t_f), \dot{x}^*(t_f), a^*, t_f) = 0$

7) $x^*(t_0) = x_0$

then $x^*(t)$ is a min-max solution of (1).

*Example 2*: Suppose that

$$J(x(t), a) = \int_0^1 \left(x^2(t) + \dot{x}^2(t) - a^2 - 2ax(t)\right) dt, \; x(0) = 0, \; x(1) = \text{free}, \; a \in \mathbb{R}$$

From condition 1 of theorem 2

$x^*(t) = Ae^t + Be^{-t} + a$.

To satisfy the boundary condition $x(0) = 0$, $A + B + a = 0$ is obtained. From condition 6

$0 = \frac{\partial g}{\partial \dot{x}}(x^*(t_f), \dot{x}^*(t_f), a^*, t_f) = 2\dot{x}^*(t_f) = Ae^{t_f} - Be^{-t_f} = Ae - B/e$,

therefore

$x^*(t) = \frac{-a}{e^2 + 1} e^t + \frac{-ae^2}{e^2 + 1} e^{-t} + a$.

On the other hand,

$$\frac{\partial g}{\partial a}(x(t),\dot{x}(t),a,t) = -2(x+a).$$

Thus, condition 2 is satisfied if

$$\int_0^1 \frac{\partial g}{\partial a}(x^*(t),\dot{x}^*(t),a^*,t)dt = 0.$$

Hence, it follows that $a^* = 0$.

Therefore, the extremum solution is

$$x^*(t) = 0$$

and

$$J(x^*,a^*) = 0.$$

Also,

$$\int_0^1 \frac{\partial^2 g}{\partial a^2}(x^*(t),\dot{x}^*(t),a^*,t)dt = -2 < 0,$$

$$\frac{\partial^2 g}{\partial \dot{x}^2}(x^*(t),\dot{x}^*(t),a^*,t) = 2 > 0$$

and

$$\frac{\partial^2 g}{\partial x^2}(x^*(t),\dot{x}^*(t),a^*,t) = 2 > 0.$$

So, all conditions are satisfied and the above extremum solution is a min-max trajectory of the functional.

### C. Free final time and specified final point

In this case, initial state and final point are specified but final time is free. Necessary conditions in this case are given in theorem 3 (proved in Appendix III).

*Theorem 3*: Consider uncertain functional (1) subject to

$$x(t_0) = x_0, x(t_f) = x_f, \; t_f > t_0 \text{ is free}, \; t \in [t_0,t_f], a \in A \tag{5}$$

Suppose $J(x(t),a)$ has a min-max solution with assumed conditions and let $x^*(t)$ be admissible. If there exists a $a^* \in A$ such that

1) $\dfrac{\partial g}{\partial x}(x^*(t), \dot{x}^*(t), a^*, t) - \dfrac{d}{dt}[\dfrac{\partial g}{\partial \dot{x}}(x^*(t), \dot{x}^*(t), a^*, t)] = 0$

2) $\displaystyle\int_{t_0}^{t_f} \dfrac{\partial g}{\partial a}(x^*(t), \dot{x}^*(t), a^*, t) dt = 0$

3) $\displaystyle\int_{t_0}^{t_f} \dfrac{\partial^2 g}{\partial a^2}(x^*(t), \dot{x}^*(t), a^*, t) dt < 0$

4) $\dfrac{\partial^2 g}{\partial x^2}(x^*(t), \dot{x}^*(t), a^*, t) > 0$

5) $\dfrac{\partial^2 g}{\partial \dot{x}^2}(x^*(t), \dot{x}^*(t), a^*, t) > 0$

6) $g(x^*(t_f), \dot{x}^*(t_f), a^*, t_f) - \dfrac{\partial g}{\partial \dot{x}}(x^*(t_f), \dot{x}^*(t_f), a^*, t_f) \dot{x}(t_f) = 0$

7) $x^*(t_0) = x_0$

then $x^*(t)$ is a min-max solution of (1).

*Example 3*: suppose that

$$J(x(t), a) = \int_0^{t_f} (6ax(t) + \dfrac{1}{2}\dot{x}^2(t) + 4 - 24t - 6a^2 t + 12at + 2a + \dfrac{9}{2}a^2 t^2) dt, x(0) = 0,$$

$x(t_f) = 0, a \in R, t_f > 0$.

Following condition 1 of theorem 3

$x^*(t) = 3a^2 t^2 + At + B$

To satisfy the boundary conditions, $x(0) = 0$, $x(t_f) = 0$, $B = 0$ and $A = -3at_f$ are obtained. From conditions 2 and 6

$-6at_f^2 + 6t_f^2 + 2t_f = 0$

and

$4 - 24t_f - 6a^2 t_f + 12at_f + 2a = 0$.

By solving these equations final time and uncertain parameter are

$$a^* = 2, t_f = \frac{1}{3}.$$

Therefore

$$.x^*(t) = 6t^2 - 2t$$

is the extremum solution of the example. By using this solution in the functional

$$J(x^*, a^*) = \frac{4}{3}.$$

Also,

$$\int_0^1 \frac{\partial^2 g}{\partial a^2}(x(t), \dot{x}(t), a, t) dt = \frac{-5}{9} < 0,$$

$$\frac{\partial^2 g}{\partial \dot{x}^2}(x(t), \dot{x}(t), a, t) = 1 > 0$$

and

$$\frac{\partial^2 g}{\partial x^2}(x(t), \dot{x}(t), a, t) = 0.$$

So, all conditions are satisfied and the above extremum solution is a min-max solution of the functional.

In the example 3, $g(x^*(t), \dot{x}^*(t), a^*, t)$ depends on first order of $x(t)$ therefore $O_1(x(t), \dot{x}(t), a, \delta x(t), \delta \dot{x}(t), \delta a, t)$ in (35) in Appendix I and $\delta x(t)$ are independent. Thus, in the example condition 4 of the theorem can be replaced by

$$\frac{\partial^2 g}{\partial x^2}(x(t), \dot{x}(t), a, t) = 0.$$

D. *Free final point ant time*

In this case suppose initial state is specified and final state and final time are free. In theorem 4, which is proved in Appendix IV, necessary conditions are presented.

*Theorem 4*: Consider uncertain functional (1) subject to

$$x(t_0) = x_0, x(t_f) \text{ and } t_f \text{ are free3, } t_f > t_0, t \in [t_0, t_f], a \in A. \tag{6}$$

Suppose $J(x(t),a)$ has a min-max solution with assumed conditions and let $x^*(t)$ be admissible. If there exists a $a^* \in A$ such that

1) $\dfrac{\partial g}{\partial x}(x^*(t), \dot{x}^*(t), a^*, t) - \dfrac{d}{dt}[\dfrac{\partial g}{\partial \dot{x}}(x^*(t), \dot{x}^*(t), a^*, t)] = 0$

2) $\displaystyle\int_{t_0}^{t_f} \dfrac{\partial g}{\partial a}(x^*(t), \dot{x}^*(t), a^*, t) dt = 0$

3) $\displaystyle\int_{t_0}^{t_f} \dfrac{\partial^2 g}{\partial a^2}(x^*(t), \dot{x}^*(t), a^*, t) dt < 0$

4) $\dfrac{\partial^2 g}{\partial x^2}(x^*(t), \dot{x}^*(t), a^*, t) > 0$

5) $\dfrac{\partial^2 g}{\partial \dot{x}^2}(x^*(t), \dot{x}^*(t), a^*, t) > 0$

6) $g(x^*(t_f), \dot{x}^*(t_f), a^*, t_f) = 0$

7) $\dfrac{\partial g}{\partial \dot{x}}(x^*(t_f), \dot{x}^*(t_f), a^*, t_f) \dot{x}(t_f) = 0$

8) $x^*(t_0) = x_0$

then $x^*(t)$ is a min-max solution of (1).

*Example 4*: suppose that

$$J(x(t), a) = \int_0^{t_f} (6ax(t) + \dfrac{1}{2}\dot{x}^2(t) + 4 - 24t - 6a^2 t + 12at + 2a + 18a^2 t^2) dt,$$

$x(0) = 0$, $x(t_f) = $ free

Following condition 1 of theorem 4

$x^*(t) = 3a^2 t^2 + At + B$.

Since $x^*(0) = 0$, $B$ is equal to zero. To satisfy condition 7 of the theorem $A = -6at_f$ is obtained.

Following conditions 2 and 6

$$4 - 24t_f - 6a^2 t_f + 12at_f + 2a = 0$$

and

$$-6at_f^2 + 6t_f^2 + 2t_f = 0.$$

By solving these equations final time and uncertain parameter are

$$a^* = 2, t_f = \frac{1}{3}.$$

Therefore

$$x^*(t) = 6t^2 - 4t.$$

So $x^*(t)$ is the extremum solution. Using this solution in the functional

$$J(x^*, a^*) = \frac{4}{3}.$$

On the other hand,

$$\int_0^1 \frac{\partial^2 g}{\partial a^2}(x(t), \dot{x}(t), a, t) dt = \frac{-2}{9} < 0,$$

$$\frac{\partial^2 g}{\partial \dot{x}^2}(x(t), \dot{x}(t), a, t) = 1 > 0$$

and

$$\frac{\partial^2 g}{\partial x^2}(x(t), \dot{x}(t), a, t) = 0.$$

Also, $g(x(t), \dot{x}(t), a, t)$ depends on first order of $x(t)$. Similar to example 3, $g(x^*(t), \dot{x}^*(t), a^*, t)$ is depend on first order of $x(t)$. So, all conditions are satisfied and the extremum solution is a min-max solution of the functional.

*Remark*: If a max-max solution of the problem is desired, conditions 3-5 in theorems 1-4 are exchanged with

3) $\int_{t_0}^{t_f} \frac{\partial^2 g}{\partial a^2}(x^*(t), \dot{x}^*(t), a^*, t) dt < 0$

4) $\dfrac{\partial^2 g}{\partial x^2}(x^*(t), \dot{x}^*(t), a^*, t) < 0$

5) $\dfrac{\partial^2 g}{\partial \dot{x}^2}(x^*(t), \dot{x}^*(t), a^*, t) < 0$

Similarly if max-min solution is concerned, conditions 3-5 of the theorems are replaced with

3) $\displaystyle\int_{t_0}^{t_f} \dfrac{\partial^2 g}{\partial a^2}(x^*(t), \dot{x}^*(t), a^*, t)\,dt > 0$

4) $\dfrac{\partial^2 g}{\partial x^2}(x^*(t), \dot{x}^*(t), a^*, t) < 0$

5) $\dfrac{\partial^2 g}{\partial \dot{x}^2}(x^*(t), \dot{x}^*(t), a^*, t) < 0$

and for min-min solution

3) $\displaystyle\int_{t_0}^{t_f} \dfrac{\partial^2 g}{\partial a^2}(x^*(t), \dot{x}^*(t), a^*, t)\,dt > 0$

4) $\dfrac{\partial^2 g}{\partial x^2}(x^*(t), \dot{x}^*(t), a^*, t) > 0$

5) $\dfrac{\partial^2 g}{\partial \dot{x}^2}(x^*(t), \dot{x}^*(t), a^*, t) > 0$.

## 4. Necessary conditions with uncertainty in state equation and cost function

Let $x(t)$ and $u(t)$ be scalar functions with continuous first derivatives. The problem is to find the control $u^*(t)$ that causes the system

$$\dot{x}(t) = f(x(t), u(t), a, t), \quad x(t_0) = x_0 \tag{7}$$

to follow trajectory $x^*(t)$ that minimizes cost function

$$J(u) = h(x(t_f), t_f) + \int_{t_0}^{t_f} g(x(t), u(t), a, t)\,dt \tag{8}$$

where $a = a^*$ maximizes the cost function (i.e. control $u^*(t)$ is desired with respect to min-max solution

$(x^*(t), a^*))$, $g: R^n \times R^m \times A \times [t_0, t_f] \to R$ is a given function. It is assumed that $g$ has continuous first and second partial derivatives with respect to all of its argument. $h: R^n \times [t_0, t_f] \to R$, is a given known differentiable function emphasizing ending point. Also, $f$ has continuous first and second partial derivatives respect to all of its argument. Initial time and state $t_0$ and $x_0$ are fixed and final time and state $t_f$ and $x_f$ are specified or free respect to problem statement. $a \in A \subset R^m$ is an uncertain parameter in which $A$ is known compact set.

Note that uncertain parameter $a$ may appear in state equation or cost function according to the problem statement. Also some terms of uncertain vector $a$ may appear in cost function in which some of them appear in state equation.

If parameter $a$ is specified then the problem reduces to the optimal control problem where it is desired to find conditions satisfying $u^*(t)$ is a minimizer [24]. However, the min–max solution of the functional can be formulated in such a way that the operation of the maximization is taken over the set of uncertainty and the operation of the minimization is taken over trajectory . Therefore, the problem is to find an admissible optimal min-max solution $(u^*(t), a^*)$ which satisfying [12]

$$\max_{a \in A} J(u^*(t), a) \leq J(u^*(t), a^*) \leq \min_{a \in A} J(u(t), a^*) \text{ for all } u(t) \in R^m \text{ and } a \in A. \tag{9}$$

Both theorems 5 and 6 follow standard method in optimal control and proves of them are similar to the same cases in [24].

In the following theorem, which is proved in Appendix V, necessary conditions to find min-max solution of the functional are presented when final state and final time is specified or free. Let

$$\begin{aligned} g_f(x(t), \dot{x}(t), p(t), a, t) = g(x(t), u(t), a, t) + p^T(t)[f(x(t), u(t), a, t) - \dot{x}(t)] + \\ [\frac{\partial h}{\partial x}(x(t), t)]^T \dot{x}(t) + \frac{\partial h}{\partial t}(x(t), t) \end{aligned} \tag{10}$$

and

$$H(x(t), u(t), p(t), a, t) = g(x(t), u(t), a, t) + p^T(t)[f(x(t), u(t), a, t)] \tag{11}$$

where $H$ is the expanded Hamiltonian function [24].

*Theorem 5*: consider the state equation (7) with cost function (8). Suppose $J(u(t),a)$ has a min-max solution with assumed conditions and let $u^*(t)$ be admissible. If there exists a $a^* \in A$ such that

1) $\dot{x}^*(t) = \dfrac{\partial H}{\partial p}(x^*(t), u^*(t), p^*(t), a^*, t)$

2) $\dot{p}^*(t) = -\dfrac{\partial H}{\partial x}(x^*(t), u^*(t), p^*(t), a^*, t)$

3) $0 = \dfrac{\partial H}{\partial u}(x^*(t), u^*(t), p^*(t), a^*, t)$

4) $0 = \dfrac{\partial H}{\partial a}(x^*(t), u^*(t), p^*(t), a^*, t)$

5) $\displaystyle\int_{t_0}^{t_f} \dfrac{\partial^2 g_f}{\partial a^2}(x^*(t), \dot{x}^*(t), a^*, t)\,dt < 0$

6) $\dfrac{\partial^2 g_f}{\partial x^2}(x^*(t), \dot{x}^*(t), a^*, t) > 0$

7) $\dfrac{\partial^2 g_f}{\partial \dot{x}^2}(x^*(t), \dot{x}^*(t), a^*, t) > 0$

8a) $x^*(t_f) = x_f$

Then $u^*(t)$ is a min-max solution of (7). If $t_f$ is fixed and $x(t_f)$ is free, then condition 8a is changed to

8b) $\dfrac{\partial h}{\partial x}(x^*(t_f), t_f) - p^*(t_f) = 0$

else if $t_f$ is free and $x(t_f)$ is fixed, condition 8a is transferred to

8c) $\dfrac{\partial h}{\partial t}(x^*(t_f), t_f) + H(x^*(t_f), u^*(t_f), p^*(t_f), a, t_f) = 0$

If both $t_f$ and $x(t_f)$ are free and independent, in order to have min-max solution, instead of condition 8a

8d) $\dfrac{\partial h}{\partial x}(x^*(t_f), t_f) - p^*(t_f) = 0$ and $\dfrac{\partial h}{\partial t}(x^*(t_f), t_f) + H(x^*(t_f), u^*(t_f), p^*(t_f), a, t_f) = 0$

must be satisfied; but if $t_f$ and $x(t_f)$ are related by

$$x(t) = \theta(t) \tag{12}$$

where $\theta(t)$ is a known function, then condition 8a is exchanged to

8e) $x(t_f) = \theta(t_f)$ and

$$[\frac{\partial h}{\partial x}(x^*(t_f),t_f) - p^*(t_f)]^T (\frac{d\theta}{dt}(t_f)) + \frac{\partial h}{\partial t}(x^*(t_f),t_f) + H(x^*(t_f),u^*(t_f),p^*(t_f),a^*,t_f) = 0$$

## 5. Necessary conditions with uncertainty in initial state

Sometimes in an optimal problem uncertain parameter may be appeared in initial state because of inexact measuring. In this section assume that initial state is not exactly fixed and has uncertainty measurement such as

$$x(t_0) = x_0 + a \tag{13}$$

where $x_0$ is known and $a \in A$ is an uncertain parameter. Suppose

$$\dot{x}(t) = f(x(t),u(t),t) \tag{14}$$

with the definition same as section 4. Using transformation equation

$$y(t) = x(t) - a \tag{15}$$

the initial state is exchanged to

$$y(t_0) = x_0 \tag{16}$$

and the cost function

$$J(u) = h_1(y(t_0),a,t_0) + \int_{t_0}^{t_f} \{g_1(y(t),u(t),a,t) + [\frac{\partial h_1}{\partial y}(x(t),a,t)]^T \dot{y}(t) + \frac{\partial h_1}{\partial t}(y(t),a,t)\}dt \tag{17}$$

where

$$h_1(y(t),a,t) = h(y(t)+a,t) \tag{18}$$

$$g_1(y(t),u(t),a,t) = g(y(t)+a,u(t),t) \tag{19}$$

and the state equation

$$\dot{y}(t) = f_1(y(t), u(t), a, t) = f(y(t) + a, u(t), t). \tag{20}$$

Let

$$g_f(y(t), \dot{y}(t), p(t), a, t) = g_1(y(t), u(t), a, t) + p^T(t)[f_1(x(t), u(t), a, t) - \dot{y}(t)] + \\ [\frac{\partial h_1}{\partial y}(y(t), a, t)]^T \dot{y}(t) + \frac{\partial h_1}{\partial t}(y(t), a, t). \tag{21}$$

In the following theorem, which is proved in Appendix VI, necessary conditions to find min-max solution of the uncertain problem where uncertainty in initial state is transferred in state equation and cost function.

*Theorem 6*: consider uncertain functional (17) with the state equation (20) subject to

$$y(t_0) = x_0, y(t_f) = y_f, t \in [t_0, t_f], a \in A. \tag{22}$$

Suppose $J(u(t), a)$ has a min-max solution with assumed conditions and let $u^*(t)$ be admissible. If there exists a $a^* \in A$ such that

1) $\dot{y}^*(t) = f_1(y^*(t), u^*(t), a^*, t)$

2) $\dot{p}^*(t) = -\frac{\partial}{\partial y}\{[f_1(y^*(t), u^*(t), a^*, t)]^T p^*(t)\} - \frac{\partial g_1}{\partial y}(y^*(t), u^*(t), a^*, t)$

3) $\frac{\partial g_1}{\partial u}(y^*(t), u^*(t), a^*, t) + \frac{\partial}{\partial u}\{[f_1(y^*(t), u^*(t), a^*, t)]^T p^*(t)\} = 0$

4) $\frac{\partial}{\partial a}\{[f_1(y^*(t), u^*(t), a^*, t)]^T p^*(t)\} + \frac{\partial}{\partial a}[\frac{\partial h_1}{\partial t}(y^*(t), a^*, t)] + \frac{\partial}{\partial a}\{[(\frac{\partial h_1}{\partial y}(y^*(t), a^*, t))]^T \dot{y}^*(t)\} + \\ \frac{\partial g_1}{\partial a}(y^*(t), u^*(t), a^*, t) + \frac{\partial h_1}{\partial a}(y(t_0), a^*, t_0) = 0$

5) $\int_{t_0}^{t_f} \frac{\partial^2 g_f}{\partial a^2}(y^*(t), \dot{y}^*(t), a^*, t)dt + \frac{\partial^2 h_1}{\partial a^2}(y(t_0), a, t_0) < 0$

6) $\frac{\partial^2 g_f}{\partial y^2}(y^*(t), \dot{y}^*(t), a^*, t) > 0$

7) $\dfrac{\partial^2 g_f}{\partial \dot{y}^2}(y^*(t), \dot{y}^*(t), a^*, t) > 0$

and boundary condition

8a) $x^*(t_f) = x_f$

then $u^*(t)$ is a min-max solution of (17). If $t_f$ is fixed and $y(t_f)$ is free, then condition 8a is changed to

8b) $\dfrac{\partial h_1}{\partial y}(y^*(t_f), a^*, t_f) - p^*(t_f) = 0$

else if $t_f$ is free and $y(t_f)$ is fixed, condition 8a is transferred to

8c) $g_1(y^*(t_f), u^*(t_f), a^*, t_f) + \dfrac{\partial h_1}{\partial t}(y^*(t_f), a^*, t_f) + p^{*T}(t_f) f_1(y^*(t_f), u^*(t_f), a^*, t_f) = 0$.

If both $t_f$ and $y(t_f)$ are free and independent, in order to have min-max solution, instead of condition 8a

8d) $\dfrac{\partial h_1}{\partial y}(y^*(t_f), a^*, t_f) - p^*(t_f) = 0$ and

$g_1(y^*(t_f), u^*(t_f), a^*, t_f) + \dfrac{\partial h_1}{\partial t}(y^*(t_f), a^*, t_f) + p^{*T}(t_f) f_1(y^*(t_f), u^*(t_f), a^*, t_f) = 0$

must be satisfied; but if $t_f$ and $y(t_f)$ are related by $y(t) = \theta(t)$ where $\theta(t)$ is a known function, then condition 8a is exchanged to

8e) $y(t_f) = \theta(t_f)$ and

$[\dfrac{\partial h_1}{\partial y}(y^*(t_f), a^*, t_f) - p^*(t_f)]^T (\dfrac{d\theta}{dt}(t_f)) + g_1(y^*(t_f), u^*(t_f), a^*, t_f) +$

$+ \dfrac{\partial h_1}{\partial t}(y^*(t_f), a^*, t_f) + p^{*T}(t_f) f_1(y^*(t_f), u^*(t_f), a^*, t_f) = 0$

*Remark*: Note that if the uncertainty parameter $a$ given in all theorems is not in the boundary defined on the problem statement (i.e. $a \notin A$), from Pontryagin's minimum principle [24], the correct solution for uncertain parameter is on its boundary.

*Example 5*: The system

$$\dot{x}_1(t) = x_2(t) + a^2 + 2a$$
$$\dot{x}_2(t) = -x_2(t) + u(t)$$

subject to

$$x(0) = [1\ 1]^T,\ x(2) = \text{free}$$

is to be controlled so that its control effort is conserved; that is, the performance measure

$$J(u) = 0.5(x_1^2(2) + x_2^2(2)) + \int_0^2 (0.5u^2(t) - 10a^4)dt$$

is to be minimized respect to $u^*(t)$ and to be maximized respect to $a^*$. The admissible states and controls are not bounded. The admissible uncertain parameter is considered in two cases. The problem is to find min-max optimal control. Note that in this example uncertainty is appeared both in state equation and its cost function.

a. *Uncertain parameter is unbounded*

The first step is to form the expanded Hamiltonian

$$H(x(t), u(t), p(t), a) = 0.5u^2(t) - 10a^4 + p_1(t)x_2(t) + p_1(t)a^2 + 2p_1(t)a - p_2(t)x_2(t) + p_2(t)u(t).$$

From conditions 2, 3 and 4 of theorem 5 necessary conditions for min-max optimal control are

$$\dot{p}_1^*(t) = -\frac{\partial H}{\partial x_1} = 0$$

$$\dot{p}_2^*(t) = -\frac{\partial H}{\partial x_2} = -p_1^*(t) + p_2^*(t)$$

$$0 = \frac{\partial H}{\partial u} = u^*(t) + p_2^*(t)$$

and

$$0 = \frac{\partial H}{\partial a} = -40a^{*3} + p_1^*(t)[2 + 2a].$$

Then $u^*(t)$ and $a^*$ is solved by

$$u^*(t) = -p_2^*(t)$$

and

$$p_1^*(t) = \frac{40a^3}{2+2a}.$$

Substituting $u^*(t)$ in state equations

$$\dot{x}_1(t) = x_2(t) + a^2 + 2a$$
$$\dot{x}_2(t) = -x_2(t) - p_2^*(t).$$

Since $t_f = 2$ is fixed and $x(t_f)$ is free, from condition 8b of theorem 5 boundary conditions are evaluated as

$$\frac{\partial h}{\partial x}(x^*(t_f), t_f) - p^*(t_f) = 0$$

and

$$x(0) = [1\ 1]^T.$$

By solving above state equation ($x_1(t)$ and $x_2(t)$) and costate equation ($p_1(t)$ and $p_2(t)$) subject to boundary conditions

$$x_1^*(t) \cong -1.429t + 4.581 + 0.039e^{t+2} - 0.034e^t - 0.074e^{-t+4} + 0.025e^{-t+2}$$
$$x_2^*(t) \cong -3.092 + 0.039e^{t+2} - 0.034e^t + 0.074e^{-t+4} - 0.024e^{-t+2}$$
$$p_1^*(t) \cong 3.092$$
$$p_2^*(t) \cong 3.092 - 0.078e^{t+2} + 0.069e^t$$

Therefore

$$a^* \cong 0.632.$$

The min-max control is

$$u^*(t) = -p_2^*(t) \cong -3.092 + 0.078e^{t+2} - 0.069e^t$$

and the corresponding cost is

$$J(u^*, a^*) \cong 4.79.$$

*b. Uncertain parameter is bounded*

Suppose uncertain parameter is bounded in $[-0.5 \ 0.5]$ (i.e. $a \in A = [-0.5 \ 0.5]$). Since $a^*$ in previous case is not in the bounded state $A = [-0.5 \ 0.5]$, then from Pontryagin's minimum principle the $a^*$ is on its boundary (i.e. $a^* = -0.5$ or $a^* = 0.5$). If $a = -0.5$, state equation and costate equation are obtained as [24]

$$\dot{x}_1(t) = x_2(t) - 0.75$$
$$\dot{x}_2(t) = -x_2(t) - p_2^*(t)$$

$$\dot{p}_1^*(t) = 0$$
$$\dot{p}_2^*(t) = -p_1^*(t) + p_2^*(t).$$

Boundary conditions are evaluated as

$$\frac{\partial h}{\partial x}(x^*(t_f), t_f) - p^*(t_f) = 0$$

and

$$x(0) = [1 \ 1]^T.$$

By solving these equations with the boundary conditions

$$x_1(t) \cong -0.948t + 2.177 + 0.004e^{t+2} - 0.016e^t - 0.026e^{-t+4} - 0.011e^{-t+2}$$
$$x_2(t) \cong -0.198 + .004e^{t+2} - 0.016e^t + 0.02e^{-t+4} + 0.011e^{-t+2}$$
$$p_1(t) \cong 0.198$$
$$p_2(t) \cong 0.198 - 0.007e^{t+2} + 0.033e^t$$

the control is

$$u_1(t) = -p_2(t) \cong -0.198 + 0.007e^{t+2} - 0.033e^t.$$

In this case the cost is evaluated as

$$J \cong -1.211.$$

Now suppose

$$a = 0.5.$$

Similar to $a = -0.5$, the solution is

$x_1(t) \cong -1.347t + 4.169 + 0.033e^{t+2} - 0.031e^t - 0.064e^{-t+4} + 0.018e^{-t+2}$

$x_2(t) \cong -2.597 + 0.033e^{t+2} - .031e^t + 0.064e^{-t+4} - 0.0184e^{-t+2}$

$p_1(t) \cong 2.597$

$p_2(t) \cong 2.597 - 0.066e^{t+2} + 0.063e^t$.

The control is

$u_2(t) = -p_2(t) \cong -2.597 + 0.066e^{t+2} - 0.063e^t$

and the corresponding cost is evaluated as

$J \cong 4.379$.

Now the state equations are solved using $u(t)$, $u_1(t)$, $u_2(t)$ and boundary conditions:

$x_1(u(t)) \cong -3.092t + 4.581 + 0.039e^{t+2} - 0.034e^t + 0.039e^{-t+2} - 4.126e^{-t} + ta^2 + 2ta$

$x_2(u(t)) \cong -3.092 + 0.039e^{t+2} - 0.034e^t - 0.039e^{-t+2} + 4.126e^{-t}$

$x_1(u_1(t)) \cong -0.198t + 0.004e^{t+2} - 0.016e^t - 0.004e^{-t+2} - 1.215e^{-t} + ta^2 + 2ta + 2.177$

$x_2(u_1(t)) \cong -0.198 + 0.004e^{t+2} - 0.016e^t - 0.004e^{-t+2} - 1.215e^{-t}$

$x_1(u_2(t)) \cong -2.597t + 4.169 + 0.033e^{t+2} - 0.031e^t - 0.033e^{-t+2} - 3.628e^{-t} + ta^2 + 2ta$

$x_2(u_2(t)) \cong -2.597 + 0.033e^{t+2} - 0.031e^t - 0.033e^{-t+2} - 3.628e^{-t}$.

Their costs with respect to uncertainty are shown in Fig. 1.(a) and enlarge in Fig. 1.(b). Based on Fig. 1, the maximum of $J(a,u_2)$ in $a \in A = [-0.5\ 0.5]$ less than maximum of $J(a,u_1)$. Therefore, $u_2$ is the global min-max solution of the example in case B. Please note that the maximum of $J(a,u)$ is less that maximum of $J(a,u_1)$ and $J(a,u_2)$ for any $a \in R$. It verifies that $u$ is min-max solution of the example where $a \in R$.

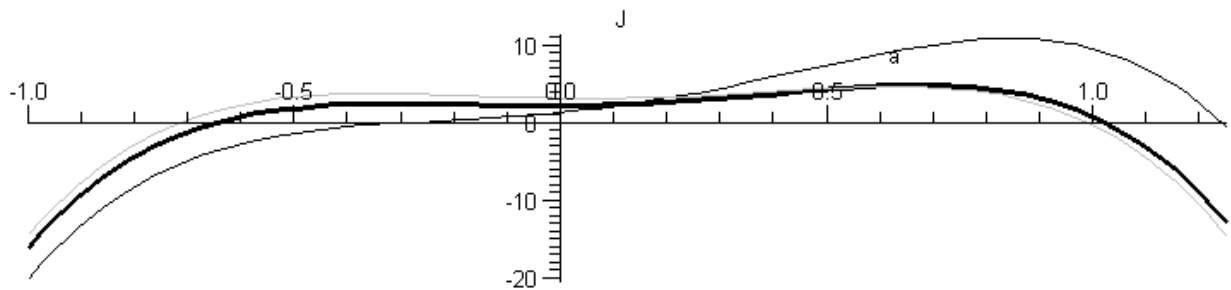

(a)

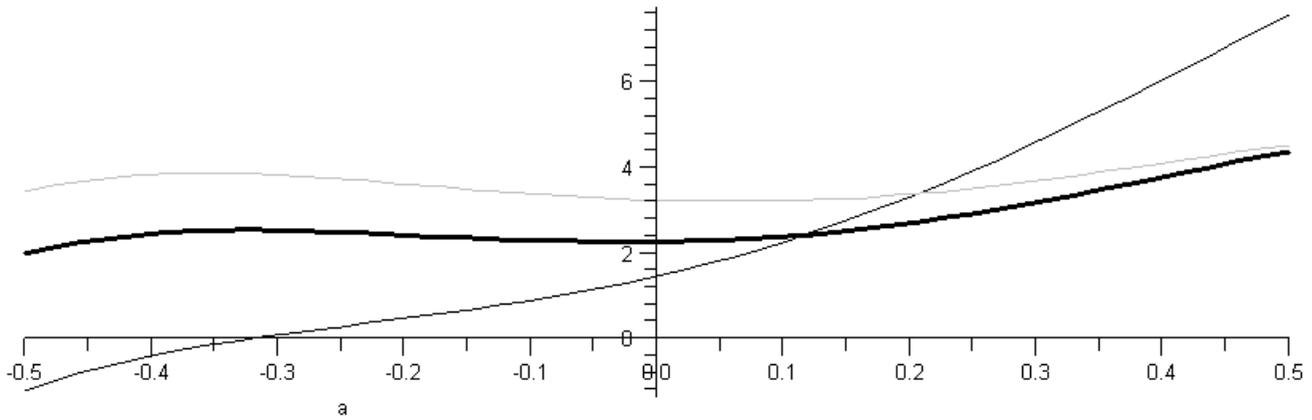

(b)

Figure 1: a. costs with respect to uncertainty, b. enlarged scale of a.

$u_1(t):$ ———, $u_2(t):$ ▬▬ $u(t):$ ———

It is clear that the min-max solution for this problem is

$$u^*(t) = -2.597 + 0.066 e^{t+2} - 0.063 e^{t}$$

where $a^* = 0.5$ and $J^* \cong 4.379$.

## 6. Conclusion

In this paper, min-max approaches of Euler-Lagrange equations for functional with uncertain parameters are formulated. Also, necessary conditions for several boundary states are developed. The problem is defined as finding necessary conditions such that the uncertain parameter maximizes the functional, where the path or trajectory is attempting to minimize. Also, the trajectory, which minimized the maximum value of the functional over all admissible uncertainty, is fined. The approach suggests necessary conditions to find min-max solution of the functional when final state and final time is specified or free.

Furthermore, in the paper the method is developed in which the functional is constraint to specific or unknown state equation. The necessary conditions to find min-max optimal control when uncertainty appears in state equation, cost function or in the initial state are presented. Finally, the method generalized when the uncertainties are bounded with using Pontryagin's minimum principle.

# Appendix I

## Proof of Theorem 1

Assume $\delta\dot{x}(t)$, $\delta x(t)$ and $\delta a$ as increments of $\dot{x}(t), x(t)$ and $a$ respectively. Note that $x(t)$ and $x(t) + \delta x(t)$ need to satisfy the boundary conditions therefore

$$\delta x(t_0) = \delta x(t_f) = 0. \tag{23}$$

Thus, the corresponding increment of the functional $J(x(t), a)$ obtains as

$$\begin{aligned}\Delta J &= J(x(t)+\delta x(t), \dot{x}(t)+\delta\dot{x}(t), a+\delta a, t) - J(x(t), \dot{x}(t), a, t) \\ &= \int_{t_0}^{t_f}[g(x(t)+\delta x(t), \dot{x}(t)+\delta\dot{x}(t), a+\delta a, t) - g(x(t), \dot{x}(t), a, t)]dt.\end{aligned} \tag{24}$$

By using Taylor's theorem for $g(x(t)+\delta x(t), \dot{x}(t)+\delta\dot{x}(t), a+\delta a, t)$

$$\begin{aligned}g(x(t)+\delta x(t), \dot{x}(t)+\delta\dot{x}(t), a+\delta a, t) &= g(x(t), \dot{x}(t), a, t) + \frac{\partial g}{\partial x}(x(t), \dot{x}(t), a, t)\delta x \\ &+ \frac{\partial g}{\partial a}(x(t), \dot{x}(t), a, t)\delta a + \frac{\partial g}{\partial \dot{x}}(x(t), \dot{x}(t), a, t)\delta\dot{x} + O(x(t), \dot{x}(t), a, \delta x(t), \delta\dot{x}(t), \delta a, t)\end{aligned} \tag{25}$$

where the $O(x(t), \dot{x}(t), a, \delta x(t), \delta\dot{x}(t), \delta a, t)$ denotes terms of order higher than 1 relative to $\delta x(t)$, $\delta\dot{x}(t)$ and $\delta a$, and other terms represent the principal linear part of the increment $\Delta J(x(t), a)$, and hence the variation of $J(x(t), a)$ is

$$\begin{aligned}\delta J(x(t), a) = \int_{t_0}^{t_f}\{[\frac{\partial g}{\partial x}(x(t), \dot{x}(t), a, t)]\delta x(t) \\ + [\frac{\partial g}{\partial \dot{x}}(x(t), \dot{x}(t), t)]\delta\dot{x}(t) + [\frac{\partial g}{\partial a}(x(t), \dot{x}(t), a, t)]\delta a\}dt.\end{aligned} \tag{26}$$

Note that $\delta x(t)$ and $\delta\dot{x}(t)$ are related by

$$\delta x(t) = \int_{t_0}^{t_f}\delta\dot{x}(t)dt + \delta x(t_0) \tag{27}$$

So by using (23) and (27) with condition 6 of the theorem 1 and chain rule for second term of (26)

$$\int_{t_0}^{t_f} \frac{\partial g}{\partial \dot{x}}(x(t), \dot{x}(t), a, t) \delta \dot{x}(t) dt = \int_{t_0}^{t_f} \frac{\partial g}{\partial \dot{x}}(x(t), \dot{x}(t), a, t) d\delta x(t)$$

$$= \frac{\partial g}{\partial \dot{x}}(x(t), \dot{x}(t), a, t) \delta x(t) \bigg|_{t_0}^{t_f} - \int_{t_0}^{t_f} \delta x(t) \frac{d}{dt} [\frac{\partial g}{\partial \dot{x}}(x(t), \dot{x}(t), a, t)] dt.$$ (28)

By substituting (28) in (26), $\delta J(x(t), a)$ is converted to

$$\delta J(x(t), a) = \frac{\partial g}{\partial \dot{x}}(x(t), \dot{x}(t), a, t) \delta x(t) \bigg|_{t_0}^{t_f} + \int_{t_0}^{t_f} [\frac{\partial g}{\partial a}(x(t), \dot{x}(t), a, t)] \delta a(t) dt$$

$$+ \int_{t_0}^{t_f} [\frac{\partial g}{\partial x}(x(t), \dot{x}(t), a, t) - \frac{d}{dt}(\frac{\partial g}{\partial \dot{x}}(x(t), \dot{x}(t), a, t)] \delta x(t) dt.$$ (29)

On the other hand,

$$\frac{\partial g}{\partial \dot{x}}(x(t), \dot{x}(t), a, t) \delta x(t) \bigg|_{t_0}^{t_f} = \frac{\partial g}{\partial \dot{x}}(x(t_f), \dot{x}(t_f), a, t_f) \delta x(t_f) - \frac{\partial g}{\partial \dot{x}}(x(t_0), \dot{x}(t_0), a, t_0) \delta x(t_0).$$ (30)

Hence, by substituting (23) in (30)

$$\frac{\partial g}{\partial \dot{x}}(x(t), \dot{x}(t), a, t) \delta x(t) \bigg|_{t_0}^{t_f} = 0.$$ (31)

A necessary condition to existence an extremum for functional $J(x(t), a)$ in $x^*(t)$ is that its variation vanishes for $x(t) = x^*(t)$, i.e., $\delta J(x(t), a) = 0$ for $x(t) = x^*(t)$ and all admissible $t \in [t_0, t_f]$. Therefore, by inserting (31) into (29) for extremum $x(t) = x^*(t)$,

$$\delta J(x^*(t), a^*) = \int_{t_0}^{t_f} [\frac{\partial g}{\partial a}(x^*(t), \dot{x}^*(t), a^*, t)] \delta a(t) dt$$

$$+ \int_{t_0}^{t_f} [\frac{\partial g}{\partial x}(x^*(t), \dot{x}^*(t), a^*, t) - \frac{d}{dt}(\frac{\partial g}{\partial \dot{x}}(x^*(t), \dot{x}^*(t), a^*, t)] \delta x(t) dt = 0$$ (32)

is obtained. Note that (32) is fulfilled for every admissible $\delta a$ and $\delta x(t)$ if and only if conditions 1 and 2 of theorem are satisfied. The min-max solution of (1) must satisfy (2). For weak conditions, it is necessary that

$$J(x^*(t), a^* + \delta a) \leq J(x^*(t), a^*) \leq J(x^*(t) + \delta x, a^*) \quad (33)$$

or equivalently

$$\Delta J(x^*(t), a^* + \delta a) \leq 0 \leq \Delta J(x^*(t) + \delta x, a^*). \tag{34}$$

Using Taylor's Theorem with some second high order term as

$$\begin{aligned}
\Delta J(x(t), a) &= \delta J(x(t), a) + [\frac{\partial^2 g}{\partial \dot{x}^2}(x(t), \dot{x}(t), t)](\delta \dot{x}(t))^2 \\
&\quad + [\frac{\partial g^2}{\partial a^2}(x(t), \dot{x}(t), a, t)](\delta a)^2 \} dt + O_1(x(t), \dot{x}(t), a, \delta x(t), \delta \dot{x}(t), \delta a, t) \\
&\cong \int_{t_0}^{t_f} \{[\frac{\partial^2 g}{\partial x^2}(x(t), \dot{x}(t), a, t)](\delta x(t))^2 + [\frac{\partial^2 g}{\partial \dot{x}^2}(x(t), \dot{x}(t), t)](\delta \dot{x}(t))^2 \} dt \\
&\quad + (\delta a)^2 \int_{t_0}^{t_f} [\frac{\partial g^2}{\partial a^2}(x(t), \dot{x}(t), a, t)] dt.
\end{aligned} \tag{35}$$

Note that $\delta a$ and $t$ are independent. Thus, based on conditions 3, 4 and 5 of the theorem (35) is hold.

∎

**Appendix II**

**Proof of Theorem2**

Following the proof in case A, it is straightforward to verify that conditions 1- 5 and 7 are given. Similar to (29), the variation of $J(x(t), a)$ is

$$\begin{aligned}
\delta J(x(t), a) &= \frac{\partial g}{\partial \dot{x}}(x^*(t_f), \dot{x}^*(t_f), a^*, t_f) \delta x(t_f) - \frac{\partial g}{\partial \dot{x}}(x^*(t_0), \dot{x}^*(t_0), a^*, t_0) \delta x(t_0) \\
&\quad + \int_{t_0}^{t_f} [\frac{\partial g}{\partial x}(x(t), \dot{x}(t), a, t) - \frac{d}{dt}(\frac{\partial g}{\partial \dot{x}}(x(t), \dot{x}(t), a, t)] \delta x(t) dt + \int_{t_0}^{t_f} [\frac{\partial g}{\partial a}(x(t), \dot{x}(t), a, t)] \delta a(t) dt = 0
\end{aligned} \tag{36}$$

Conditions 1, 2 and 7 imply all the terms in the right-hand side of the equation (36) to zero except the first term. Therefore to have an extremum path condition 6 must be satisfied. Also, conditions 3, 4 and 5 are given from (35).

∎

## Appendix III

## Proof of theorem 3

The corresponding increment $\Delta J(x(t),a)$ of the cost $J(x(t),a)$ is (Kirk, 1970)

$$\Delta J = \int_{t_0}^{t_f+\delta t_f} g(x(t),\dot{x}(t),a,t)dt - \int_{t_0}^{t_f} g(x^*(t),\dot{x}^*(t),a^*,t)dt =$$

$$\int_{t_0}^{t_f} \{g(x^*(t)+\delta x(t),\dot{x}^*(t)+\delta\dot{x}(t),a^*+\delta a,t) - g(x^*(t),\dot{x}^*(t),a^*,t)\}dt + \int_{t_f}^{t_f+\delta t_f} g(x(t),\dot{x}(t),a,t)dt. \quad (37)$$

It follows by using Taylor's theorem,

$$\Delta J = \int_{t_0}^{t_f} \{[\frac{\partial g}{\partial x}(x^*(t),\dot{x}^*(t),a^*,t)]\delta x(t) + [\frac{\partial g}{\partial \dot{x}}(x^*(t),\dot{x}^*(t),a^*,t)]\delta\dot{x}(t)$$

$$+ \frac{\partial g}{\partial a}(x^*(t),\dot{x}^*(t),a^*,t)]\delta a\}dt + O(\cdot) + \int_{t_f}^{t_f+\delta t_f} g(x(t),\dot{x}(t),a,t)dt \quad (38)$$

where the $O(\cdot)$ denotes terms of order higher than 1 relative to $\delta x(t)$, $\delta\dot{x}(t)$ and $\delta a$. The second integral is rewrite as

$$\int_{t_f}^{t_f+\delta t_f} g(x(t),\dot{x}(t),a,t)dt = [g(x(t_f),\dot{x}(t_f),a,t_f)]\delta t_f + O'(\cdot) \quad (39)$$

in which $O'(\cdot)$ denotes terms of order higher than 1 relative to $\delta t_f$. It is clear that

$$\delta x(t_f) = -\dot{x}^*(t_f)\delta t_f \quad (40)$$

therefore

$$0 = \delta J(x^*(t),a^*) = \{[-\frac{\partial g}{\partial \dot{x}}(x^*(t_f),\dot{x}^*(t_f),a^*,t_f)]\dot{x}^*(t_f) + g(x^*(t_f),\dot{x}^*(t_f),a^*,t_f)\}\delta t_f$$

$$+ \int_{t_0}^{t_f} \{[\frac{\partial g}{\partial x}(x^*(t),\dot{x}^*(t),a^*,t) - \frac{d}{dt}[\frac{\partial g}{\partial \dot{x}}(x^*(t),\dot{x}^*(t),a^*,t)]]\delta x(t)\}dt + \int_{t_0}^{t_f} \frac{\partial g}{\partial a}(x^*(t),\dot{x}^*(t),a^*,t)\delta a dt. \quad (41)$$

In order to verify equation (41), conditions 1, 2, 6 and 7 must be satisfied. Similar to theorem 1, conditions 3, 4 and 5 must be satisfied to have a min-max solution based on (35).

∎

## Appendix IV

## Proof of Theorem 4

Following the proof in case C and the fact that

$$\delta x(t_f) = \delta x_f - \dot{x}^*(t_f)\delta t_f. \tag{42}$$

Then, by substituting (36) and (42) in (38) we have

$$\begin{aligned}
0 = \delta J(x^*, a^*) &= \frac{\partial g}{\partial \dot{x}}(x^*(t_f), \dot{x}^*(t_f), a^*, t_f)\delta x_f + [g(x^*(t_f), \dot{x}^*(t_f), a^*, t_f) \\
&- \frac{\partial g}{\partial \dot{x}}(x^*(t_f), \dot{x}^*(t_f), a^*, t_f)\dot{x}^*(t_f)]\delta t_f + \int_{t_0}^{t_f}\{[\frac{\partial g}{\partial x}(x^*(t), \dot{x}^*(t), a^*, t) \\
&- \frac{d}{dt}[\frac{\partial g}{\partial \dot{x}}(x^*(t), \dot{x}^*(t), a^*, t)]]\delta x(t)\}dt + \int_{t_0}^{t_f}\frac{\partial g}{\partial a}(x^*(t), \dot{x}^*(t), a^*, t)\delta a dt.
\end{aligned} \tag{43}$$

Since $t_f$ and $x(t_f)$ are independent therefore

$$\frac{\partial g}{\partial \dot{x}}(x^*(t_f), \dot{x}^*(t_f), a^*, t_f) = 0 \tag{44}$$

and

$$g(x^*(t_f), \dot{x}^*(t_f), a^*, t_f) - [\frac{\partial g}{\partial \dot{x}}(x^*(t_f), \dot{x}^*(t_f), a^*, t_f)]\dot{x}^*(t_f) = 0. \tag{45}$$

Thus

$$g(x^*(t_f), \dot{x}^*(t_f), a^*, t_f) = 0 \tag{46}$$

in which conditions 6 and 7 are verified and by following the proof of Theorem 3 other conditions will be proved. Similar to Theorem 1, conditions 3, 4 and 5 must be satisfied to have a min-max solution.

∎

## Appendix V

## Proof of Theorem 5

Since $h$ is differentiable

$$h(x(t_f), t_f) = \int_{t_0}^{t_f} \frac{d}{dt}[h(x(t),t)]dt + h(x(t_0), t_0) \tag{47}$$

So that the quadratic cost functional (8) can be expressed as

$$J(u) = h(x(t_0), t_0) + \int_{t_0}^{t_f} \{g(x(t), u(t), a, t) + \frac{d}{dt}[h(x(t),t)]\}dt. \tag{48}$$

$x(t_0)$ and $t_0$ are fixed, thus the minimizing does not affect the first term of $J$. Similar to Kirk (1970) page 185, using chain rule of differentiation, the problem is minimizing the quadratic cost function

$$J(u) = \int_{t_0}^{t_f} \{g(x(t), u(t), a, t) + [\frac{\partial h}{\partial x}(x(t), t)]^T \dot{x}(t) + \frac{\partial h}{\partial t}(x(t), t)\}dt. \tag{49}$$

Including the differential equation constraints (7) by using the Lagrange multipliers $p^T(t) = [p_1(t), ..., p_n(t)]$ the cost function is exchanged to

$$J_f(u) = \int_{t_0}^{t_f} \{g(x(t), u(t), t) + [\frac{\partial h}{\partial x}(x(t), t)]^T \dot{x}(t) + \frac{\partial h}{\partial t}(x(t), t) + p^T(t)[f(x(t), u(t), a, t) - \dot{x}(t)]\}dt. \tag{50}$$

The variation of $J_f$ is

$$\delta J_f(u^*, a^*) = \left[\frac{\partial g_f}{\partial \dot{x}}^T (x^*(t_f), \dot{x}^*(t_f), u^*(t_f), p^*(t_f), a^*, t_f)\right]\delta x_f + \{g_f(x^*(t_f), \dot{x}^*(t_f), u^*(t_f), p^*(t_f), a^*, t_f) -$$

$$\left[\frac{\partial g_f}{\partial \dot{x}}^T (x^*(t_f), \dot{x}^*(t_f), u^*(t_f), p^*(t_f), a^*, t_f)\right]\dot{x}^*(t_f)\}\delta t_f + \int_{t_0}^{t_f}\left\{\left[\frac{\partial g_f}{\partial x}^T (x^*(t_f), \dot{x}^*(t_f), u^*(t_f), p^*(t_f), a^*, t_f) - \right.\right.$$

$$\frac{d}{dt}\left(\frac{\partial g_f}{\partial \dot{x}}^T (x^*(t_f), \dot{x}^*(t_f), u^*(t_f), p^*(t_f), a^*, t_f)\right)\right]\delta x(t) + \frac{\partial g_f}{\partial a}^T (x^*(t_f), \dot{x}^*(t_f), u^*(t_f), p^*(t_f), a^*, t_f)\delta a +$$

$$\frac{\partial g_f}{\partial p}^T (x^*(t_f), \dot{x}^*(t_f), u^*(t_f), p^*(t_f), a^*, t_f)\delta p(t) + \frac{\partial g_f}{\partial u}^T (x^*(t_f), \dot{x}^*(t_f), u^*(t_f), p^*(t_f), a^*, t_f)\delta u(t)\bigg\} = 0 \tag{51}$$

where

$$g_f(x(t), \dot{x}(t), p(t), a, t) = g(x(t), u(t), a, t) + p^T(t)[f(x(t), u(t), a, t) - \dot{x}(t)] +$$
$$[\frac{\partial h}{\partial x}(x(t), t)]^T \dot{x}(t) + \frac{\partial h}{\partial t}(x(t), t). \tag{52}$$

It is clear that those terms inside the integral which involve the function $h$ are equal to zero (Kirk, 1970,

page 186), therefore

$$0 = \int_{t_0}^{t_f} \{ \frac{\partial g}{\partial x}^T (x^*(t), u^*(t), a^*, t) + p^{*T}(t)[\frac{\partial f}{\partial x}(x^*(t), u^*(t), a^*, t) - \frac{d}{dt}(p^{*T}(t))]\delta x(t) +$$

$$[\frac{\partial g}{\partial u}^T (x^*(t), u^*(t), a^*, t) + p^{*T}(t)(\frac{\partial f}{\partial u}(x^*(t), u^*(t), a^*, t))]\delta u(t) + \quad (53)$$

$$[f(x^*(t), u^*(t), a^*, t)^T]\delta p(t) + p^{*T}(t)[\frac{\partial f}{\partial a}(x^*(t), u^*(t), a^*, t)]\delta a(t) \} dt$$

This integral must vanishes on a minimum, thus

$$\dot{p}^*(t) = -[\frac{\partial f}{\partial x}(x^*(t), u^*(t), a^*, t)]^T p^*(t) - \frac{\partial g}{\partial x}(x^*(t), u^*(t), a^*, t) \quad (54)$$

$$\frac{\partial g}{\partial u}(x^*(t), u^*(t), a^*, t) + [\frac{\partial f}{\partial u}(x^*(t), u^*(t), a^*, t)]^T p^*(t) = 0 \quad (55)$$

$$[\frac{\partial f}{\partial a}(x^*(t), u^*(t), a^*, t)]^T p^*(t) = 0. \quad (56)$$

Since the variation $\delta J_t$ must be zero, the terms outside the integral satisfy

$$[\frac{\partial h}{\partial x}(x^*(t_f), t_f) - p^*(t_f)]^T \delta x_f + \{g(x^*(t_f), u^*(t_f), a^*, t_f) + \frac{\partial h}{\partial t}(x^*(t_f), t_f) +$$

$$p^*(t_f)[f(x^*(t_f), u^*(t_f), a^*, t_f)]\}^T \delta t_f = 0.$$

$$(57)$$

Using expanded Hamiltonian function defining as

$$H(x(t), u(t), p(t), a, t) = g(x(t), u(t), a, t) + p^T(t)[f(x(t), u(t), a, t)]. \quad (58)$$

The necessary conditions (7) and (54)-(56) is rewritten by following

$$\dot{x}^*(t) = \frac{\partial H}{\partial p}(x^*(t), u^*(t), p^*(t), a^*, t) \quad (59)$$

$$\dot{p}^*(t) = -\frac{\partial H}{\partial x}(x^*(t), u^*(t), p^*(t), a^*, t) \quad (60)$$

$$0 = \frac{\partial H}{\partial u}(x^*(t), u^*(t), p^*(t), a^*, t) \quad (61)$$

$$0 = \frac{\partial H}{\partial a}(x^*(t), u^*(t), p^*(t), a^*, t) \quad (62)$$

and bounded condition

$$[\frac{\partial h}{\partial x}(x^*(t_f),t_f) - p^*(t_f)]^T \delta x_f + [\frac{\partial h}{\partial t}(x^*(t_f),t_f) + H(x^*(t_f),u^*(t_f),p^*(t_f),a^*,t_f)]\delta t_f = 0. \tag{63}$$

It is clear that if final state is fixed, $\delta t_f$ and $\delta x_f$ are equal to zero, then (63) is vanished and to verify bounded condition $x(t_f) = x_f$ must be satisfied. On the other hand if $t_f$ is fixed and $x(t_f)$ is free, $\delta t_f$ is equal to zero but $\delta x_f$ is arbitrary so (63) is exchanged to

$$\frac{\partial h}{\partial x}(x^*(t_f),t_f) - p^*(t_f) = 0. \tag{64}$$

Suppose $t_f$ is free and $x(t_f)$ is fixed, then $\delta t_f$ is arbitrary but $\delta x_f$ is equal to zero. Therefore (63) is substituted by

$$\frac{\partial h}{\partial t}(x^*(t_f),t_f) + H(x^*(t_f),u^*(t_f),p^*(t_f),a,t_f) = 0. \tag{65}$$

If $t_f$ and $x(t_f)$ are free and independent, to satisfy (63)

$$\frac{\partial h}{\partial x}(x^*(t_f),t_f) - p^*(t_f) = 0 \tag{66}$$

and

$$\frac{\partial h}{\partial t}(x^*(t_f),t_f) + H(x^*(t_f),u^*(t_f),p^*(t_f),a,t_f) = 0 \tag{67}$$

must be satisfied. Now let $t_f$ and $x(t_f)$ are free and related by

$$x(t) = \theta(t) \tag{68}$$

where $\theta(t)$ is a known function, then

$$\delta x_f = [\frac{d\theta}{dt}(t_f)]\delta t_f. \tag{69}$$

By substituting (69) in (63)

$$[\frac{\partial h}{\partial x}(x^*(t_f),t_f) - p^*(t_f)]^T (\frac{d\theta}{dt}(t_f)) + \frac{\partial h}{\partial t}(x^*(t_f),t_f) + H(x^*(t_f),u^*(t_f),p^*(t_f),a^*,t_f) = 0 \tag{70}$$

## Appendix VI

## Proof of Theorem 6

Using transformation equation (15), the cost function is

$$J(u) = H(y(t_0),a,t_0) + \int_{t_0}^{t_f} \{G(y(t),u(t),a,t) + [\frac{\partial H}{\partial y}(y(t),a,t)]^T \dot{y}(t) + \frac{\partial H}{\partial t}(y(t),a,t)\}dt. \quad (71)$$

Note that the first term of $J(u)$ depends on uncertain parameter $a$ and has no efficiency on minimizing $J(u)$, but in order to have min-max solution this term appears in $\delta a$.

Following proves of theorems 5 and 6 conditions 1-8 of theorem 7 are satisfied.

## Appendix I

## Proof of Theorem 1

Assume $\delta\dot{x}(t)$, $\delta x(t)$ and $\delta a$ as increments of $\dot{x}(t), x(t)$ and $a$ respectively. Note that $x(t)$ and $x(t) + \delta x(t)$ need to satisfy the boundary conditions therefore

$$\delta x(t_0) = \delta x(t_f) = 0. \quad (23)$$

Thus, the corresponding increment of the functional $J(x(t),a)$ obtains as

$$\Delta J = J(x(t) + \delta x(t), \dot{x}(t) + \delta\dot{x}(t), a + \delta a, t) - J(x(t), \dot{x}(t), a, t)$$
$$= \int_{t_0}^{t_f} [g(x(t) + \delta x(t), \dot{x}(t) + \delta\dot{x}(t), a + \delta a, t) - g(x(t), \dot{x}(t), a, t)]dt. \quad (24)$$

By using Taylor's theorem for $g(x(t) + \delta x(t), \dot{x}(t) + \delta\dot{x}(t), a + \delta a, t)$

$$g(x(t)+\delta x(t),\dot{x}(t)+\delta\dot{x}(t),a+\delta a,t) = g(x(t),\dot{x}(t),a,t) + \frac{\partial g}{\partial x}(x(t),\dot{x}(t),a,t)\delta x$$
$$+ \frac{\partial g}{\partial a}(x(t),\dot{x}(t),a,t)\delta a + \frac{\partial g}{\partial \dot{x}}(x(t),\dot{x}(t),a,t)\delta\dot{x} + O(x(t),\dot{x}(t),a,\delta x(t),\delta\dot{x}(t),\delta a,t)$$
(25)

where the $O(x(t),\dot{x}(t),a,\delta x(t),\delta\dot{x}(t),\delta a,t)$ denotes terms of order higher than 1 relative to $\delta x(t)$, $\delta\dot{x}(t)$ and $\delta a$, and other terms represent the principal linear part of the increment $\Delta J(x(t),a)$, and hence the variation of $J(x(t),a)$ is

$$\delta J(x(t),a) = \int_{t_0}^{t_f} \{[\frac{\partial g}{\partial x}(x(t),\dot{x}(t),a,t)]\delta x(t)$$
$$+ [\frac{\partial g}{\partial \dot{x}}(x(t),\dot{x}(t),t)]\delta\dot{x}(t) + [\frac{\partial g}{\partial a}(x(t),\dot{x}(t),a,t)]\delta a\}dt.$$
(26)

Note that $\delta x(t)$ and $\delta\dot{x}(t)$ are related by

$$\delta x(t) = \int_{t_0}^{t_f} \delta\dot{x}(t)dt + \delta x(t_0)$$
(27)

So by using (23) and (27) with condition 6 of the theorem 1 and chain rule for second term of (26)

$$\int_{t_0}^{t_f} \frac{\partial g}{\partial \dot{x}}(x(t),\dot{x}(t),a,t)\delta\dot{x}(t)dt = \int_{t_0}^{t_f} \frac{\partial g}{\partial \dot{x}}(x(t),\dot{x}(t),a,t)d\delta x(t)$$
$$= \frac{\partial g}{\partial \dot{x}}(x(t),\dot{x}(t),a,t)\delta x(t)\Big|_{t_0}^{t_f} - \int_{t_0}^{t_f} \delta x(t)\frac{d}{dt}[\frac{\partial g}{\partial \dot{x}}(x(t),\dot{x}(t),a,t)]dt.$$
(28)

By substituting (28) in (26), $\delta J(x(t),a)$ is converted to

$$\delta J(x(t),a) = \frac{\partial g}{\partial \dot{x}}(x(t),\dot{x}(t),a,t)\delta x(t)\Big|_{t_0}^{t_f} + \int_{t_0}^{t_f} [\frac{\partial g}{\partial a}(x(t),\dot{x}(t),a,t)]\delta a(t)dt$$
$$+ \int_{t_0}^{t_f} [\frac{\partial g}{\partial x}(x(t),\dot{x}(t),a,t) - \frac{d}{dt}(\frac{\partial g}{\partial \dot{x}}(x(t),\dot{x}(t),a,t)]\delta x(t)dt.$$
(29)

On the other hand,

$$\frac{\partial g}{\partial \dot{x}}(x(t),\dot{x}(t),a,t)\delta x(t)\Big|_{t_0}^{t_f} = \frac{\partial g}{\partial \dot{x}}(x(t_f),\dot{x}(t_f),a,t_f)\delta x(t_f) - \frac{\partial g}{\partial \dot{x}}(x(t_0),\dot{x}(t_0),a,t_0)\delta x(t_0).$$
(30)

Hence, by substituting (23) in (30)

$$\frac{\partial g}{\partial \dot{x}}(x(t),\dot{x}(t),a,t)\delta x(t)\bigg|_{t_0}^{t_f} = 0. \tag{31}$$

A necessary condition to existence an extremum for functional $J(x(t),a)$ in $x^*(t)$ is that its variation vanishes for $x(t) = x^*(t)$, i.e., $\delta J(x(t),a) = 0$ for $x(t) = x^*(t)$ and all admissible $t \in [t_0, t_f]$. Therefore, by inserting (31) into (29) for extremum $x(t) = x^*(t)$,

$$\delta J(x^*(t), a^*) = \int_{t_0}^{t_f} [\frac{\partial g}{\partial a}(x^*(t), \dot{x}^*(t), a^*, t)]\delta a(t)dt$$
$$+ \int_{t_0}^{t_f} [\frac{\partial g}{\partial x}(x^*(t), \dot{x}^*(t), a^*, t) - \frac{d}{dt}(\frac{\partial g}{\partial \dot{x}}(x^*(t), \dot{x}^*(t), a^*, t)]\delta x(t)dt = 0 \tag{32}$$

is obtained. Note that (32) is fulfilled for every admissible $\delta a$ and $\delta x(t)$ if and only if conditions 1 and 2 of theorem are satisfied. The min-max solution of (1) must satisfy (2). For weak conditions, it is necessary that

$$J(x^*(t), a^* + \delta a) \leq J(x^*(t), a^*) \leq J(x^*(t) + \delta x, a^*) \tag{33}$$

or equivalently

$$\Delta J(x^*(t), a^* + \delta a) \leq 0 \leq \Delta J(x^*(t) + \delta x, a^*). \tag{34}$$

Using Taylor's Theorem with some second high order term as

$$\Delta J(x(t),a) = \delta J(x(t),a) + [\frac{\partial^2 g}{\partial \dot{x}^2}(x(t), \dot{x}(t), t)](\delta \dot{x}(t))^2$$
$$+ [\frac{\partial g^2}{\partial a^2}(x(t), \dot{x}(t), a, t)](\delta a)^2\}dt + O_1(x(t), \dot{x}(t), a, \delta x(t), \delta \dot{x}(t), \delta a, t)$$
$$\cong \int_{t_0}^{t_f} \{[\frac{\partial^2 g}{\partial x^2}(x(t), \dot{x}(t), a, t)](\delta x(t))^2 + [\frac{\partial^2 g}{\partial \dot{x}^2}(x(t), \dot{x}(t), t)](\delta \dot{x}(t))^2\}dt \tag{35}$$
$$+ (\delta a)^2 \int_{t_0}^{t_f} [\frac{\partial g^2}{\partial a^2}(x(t), \dot{x}(t), a, t)]dt.$$

Note that $\delta a$ and $t$ are independent. Thus, based on conditions 3, 4 and 5 of the theorem (35) is hold.

∎

## Appendix II

### Proof of Theorem2

Following the proof in case A, it is straightforward to verify that conditions 1- 5 and 7 are given. Similar to (29), the variation of $J(x(t),a)$ is

$$\delta J(x(t),a) = \frac{\partial g}{\partial \dot{x}}(x^*(t_f),\dot{x}^*(t_f),a^*,t_f)\delta x(t_f) - \frac{\partial g}{\partial \dot{x}}(x^*(t_0),\dot{x}^*(t_0),a^*,t_0)\delta x(t_0)$$
$$+ \int_{t_0}^{t_f}[\frac{\partial g}{\partial x}(x(t),\dot{x}(t),a,t) - \frac{d}{dt}(\frac{\partial g}{\partial \dot{x}}(x(t),\dot{x}(t),a,t)]\delta x(t)dt + \int_{t_0}^{t_f}[\frac{\partial g}{\partial a}(x(t),\dot{x}(t),a,t)]\delta a(t)dt = 0 \quad (36)$$

Conditions 1, 2 and 7 imply all the terms in the right-hand side of the equation (36) to zero except the first term. Therefore to have an extremum path condition 6 must be satisfied. Also, conditions 3, 4 and 5 are given from (35).

∎

## Appendix III

### Proof of theorem 3

The corresponding increment $\Delta J(x(t),a)$ of the cost $J(x(t),a)$ is (Kirk, 1970)

$$\Delta J = \int_{t_0}^{t_f+\delta t_f} g(x(t),\dot{x}(t),a,t)dt - \int_{t_0}^{t_f} g(x^*(t),\dot{x}^*(t),a^*,t)dt =$$
$$\int_{t_0}^{t_f}\{g(x^*(t)+\delta x(t),\dot{x}^*(t)+\delta \dot{x}(t),a^*+\delta a,t) - g(x^*(t),\dot{x}^*(t),a^*,t)\}dt + \int_{t_f}^{t_f+\delta t_f} g(x(t),\dot{x}(t),a,t)dt. \quad (37)$$

It follows by using Taylor's theorem,

$$\Delta J = \int_{t_0}^{t_f}\{[\frac{\partial g}{\partial x}(x^*(t),\dot{x}^*(t),a^*,t)]\delta x(t) + [\frac{\partial g}{\partial \dot{x}}(x^*(t),\dot{x}^*(t),a^*,t)]\delta \dot{x}(t)$$
$$+ \frac{\partial g}{\partial a}(x^*(t),\dot{x}^*(t),a^*,t)]\delta a\}dt + O(\cdot) + \int_{t_f}^{t_f+\delta t_f} g(x(t),\dot{x}(t),a,t)dt \quad (38)$$

where the $O(\cdot)$ denotes terms of order higher than 1 relative to $\delta x(t)$, $\delta \dot{x}(t)$ and $\delta a$. The second integral is

rewrite as

$$\int_{t_f}^{t_f + \delta t_f} g(x(t), \dot{x}(t), a, t) dt = [g(x(t_f), \dot{x}(t_f), a, t_f)] \delta t_f + O'(\cdot) \tag{39}$$

in which $O'(\cdot)$ denotes terms of order higher than 1 relative to $\delta t_f$. It is clear that

$$\delta x(t_f) = -\dot{x}^*(t_f) \delta t_f \tag{40}$$

therefore

$$0 = \delta J(x^*(t), a^*) = \{[-\frac{\partial g}{\partial \dot{x}}(x^*(t_f), \dot{x}^*(t_f), a^*, t_f)] \dot{x}^*(t_f) + g(x^*(t_f), \dot{x}^*(t_f), a^*, t_f)\} \delta t_f$$
$$+ \int_{t_0}^{t_f} \{[\frac{\partial g}{\partial x}(x^*(t), \dot{x}^*(t), a^*, t) - \frac{d}{dt}[\frac{\partial g}{\partial \dot{x}}(x^*(t), \dot{x}^*(t), a^*, t)]] \delta x(t)\} dt + \int_{t_0}^{t_f} \frac{\partial g}{\partial a}(x^*(t), \dot{x}^*(t), a^*, t) \delta a \, dt. \tag{41}$$

In order to verify equation (41), conditions 1, 2, 6 and 7 must be satisfied. Similar to theorem 1, conditions 3, 4 and 5 must be satisfied to have a min-max solution based on (35).

∎

**Appendix IV**

**Proof of Theorem 4**

Following the proof in case C and the fact that

$$\delta x(t_f) = \delta x_f - \dot{x}^*(t_f) \delta t_f. \tag{42}$$

Then, by substituting (36) and (42) in (38) we have

$$0 = \delta J(x^*, a^*) = \frac{\partial g}{\partial \dot{x}}(x^*(t_f), \dot{x}^*(t_f), a^*, t_f) \delta x_f + [g(x^*(t_f), \dot{x}^*(t_f), a^*, t_f)$$
$$- \frac{\partial g}{\partial \dot{x}}(x^*(t_f), \dot{x}^*(t_f), a^*, t_f) \dot{x}^*(t_f)] \delta t_f + \int_{t_0}^{t_f} \{[\frac{\partial g}{\partial x}(x^*(t), \dot{x}^*(t), a^*, t) \tag{43}$$
$$- \frac{d}{dt}[\frac{\partial g}{\partial \dot{x}}(x^*(t), \dot{x}^*(t), a^*, t)]] \delta x(t)\} dt + \int_{t_0}^{t_f} \frac{\partial g}{\partial a}(x^*(t), \dot{x}^*(t), a^*, t) \delta a \, dt.$$

Since $t_f$ and $x(t_f)$ are independent therefore

$$\frac{\partial g}{\partial \dot{x}}(x^*(t_f), \dot{x}^*(t_f), a^*, t_f) = 0 \tag{44}$$

and

$$g(x^*(t_f), \dot{x}^*(t_f), a^*, t_f) - [\frac{\partial g}{\partial \dot{x}}(x^*(t_f), \dot{x}^*(t_f), a^*, t_f)]\dot{x}^*(t_f) = 0. \tag{45}$$

Thus

$$g(x^*(t_f), \dot{x}^*(t_f), a^*, t_f) = 0 \tag{46}$$

in which conditions 6 and 7 are verified and by following the proof of Theorem 3 other conditions will be proved. Similar to Theorem 1, conditions 3, 4 and 5 must be satisfied to have a min-max solution.

■

## Appendix V

## Proof of Theorem 5

Since $h$ is differentiable

$$h(x(t_f), t_f) = \int_{t_0}^{t_f} \frac{d}{dt}[h(x(t), t)]dt + h(x(t_0), t_0) \tag{47}$$

So that the quadratic cost functional (8) can be expressed as

$$J(u) = h(x(t_0), t_0) + \int_{t_0}^{t_f} \{g(x(t), u(t), a, t) + \frac{d}{dt}[h(x(t), t)]\}dt. \tag{48}$$

$x(t_0)$ and $t_0$ are fixed, thus the minimizing does not affect the first term of $J$. Similar to Kirk (1970) page 185, using chain rule of differentiation, the problem is minimizing the quadratic cost function

$$J(u) = \int_{t_0}^{t_f} \{g(x(t), u(t), a, t) + [\frac{\partial h}{\partial x}(x(t), t)]^T \dot{x}(t) + \frac{\partial h}{\partial t}(x(t), t)\}dt. \tag{49}$$

Including the differential equation constraints (7) by using the Lagrange multipliers $p^T(t) = [p_1(t), ..., p_n(t)]$ the cost function is exchanged to

$$J_f(u) = \int_{t_0}^{t_f} \{g(x(t),u(t),t) + [\frac{\partial h}{\partial x}(x(t),t)]^T \dot{x}(t) + \frac{\partial h}{\partial t}(x(t),t) + p^T(t)[f(x(t),u(t),a,t) - \dot{x}(t)]\}dt. \quad (50)$$

The variation of $J_f$ is

$$\delta J_f(u^*,a^*) = \left[\frac{\partial g_f}{\partial \dot{x}}^T (x^*(t_f),\dot{x}^*(t_f),u^*(t_f),p^*(t_f),a^*,t_f)\right]\delta x_f + \{g_f(x^*(t_f),\dot{x}^*(t_f),u^*(t_f),p^*(t_f),a^*,t_f) -$$

$$\left[\frac{\partial g_f}{\partial \dot{x}}^T (x^*(t_f),\dot{x}^*(t_f),u^*(t_f),p^*(t_f),a^*,t_f)\right]\dot{x}^*(t_f)\}\delta t_f + \int_{t_0}^{t_f}\{[\frac{\partial g_f}{\partial x}^T (x^*(t_f),\dot{x}^*(t_f),u^*(t_f),p^*(t_f),a^*,t_f) -$$

$$\frac{d}{dt}\left(\frac{\partial g_f}{\partial \dot{x}}^T (x^*(t_f),\dot{x}^*(t_f),u^*(t_f),p^*(t_f),a^*,t_f)\right)]\delta x(t) + \frac{\partial g_f}{\partial a}^T (x^*(t_f),\dot{x}^*(t_f),u^*(t_f),p^*(t_f),a^*,t_f)\delta a +$$

$$\frac{\partial g_f}{\partial p}^T (x^*(t_f),\dot{x}^*(t_f),u^*(t_f),p^*(t_f),a^*,t_f)\delta p(t) + \frac{\partial g_f}{\partial u}^T (x^*(t_f),\dot{x}^*(t_f),u^*(t_f),p^*(t_f),a^*,t_f)\delta u(t)\} = 0 \quad (51)$$

where

$$g_f(x(t),\dot{x}(t),p(t),a,t) = g(x(t),u(t),a,t) + p^T(t)[f(x(t),u(t),a,t) - \dot{x}(t)] + [\frac{\partial h}{\partial x}(x(t),t)]^T \dot{x}(t) + \frac{\partial h}{\partial t}(x(t),t). \quad (52)$$

It is clear that those terms inside the integral which involve the function $h$ are equal to zero (Kirk, 1970, page 186), therefore

$$0 = \int_{t_0}^{t_f}\{\frac{\partial g}{\partial x}^T (x^*(t),u^*(t),a^*,t) + p^{*T}(t)[\frac{\partial f}{\partial x}(x^*(t),u^*(t),a^*,t) - \frac{d}{dt}(p^{*T}(t))]\delta x(t) +$$

$$[\frac{\partial g}{\partial u}^T (x^*(t),u^*(t),a^*,t) + p^{*T}(t)(\frac{\partial f}{\partial u}(x^*(t),u^*(t),a^*,t))]\delta u(t) + \quad (53)$$

$$[f(x^*(t),u^*(t),a^*,t)^T]\delta p(t) + p^{*T}(t)[\frac{\partial f}{\partial a}(x^*(t),u^*(t),a^*,t)]\delta a(t)\}dt$$

This integral must vanishes on a minimum, thus

$$\dot{p}^*(t) = -[\frac{\partial f}{\partial x}(x^*(t),u^*(t),a^*,t)]^T p^*(t) - \frac{\partial g}{\partial x}(x^*(t),u^*(t),a^*,t) \quad (54)$$

$$\frac{\partial g}{\partial u}(x^*(t),u^*(t),a^*,t) + [\frac{\partial f}{\partial u}(x^*(t),u^*(t),a^*,t)]^T p^*(t) = 0 \quad (55)$$

$$[\frac{\partial f}{\partial a}(x^*(t),u^*(t),a^*,t)]^T p^*(t) = 0. \quad (56)$$

Since the variation $\delta J_t$ must be zero, the terms outside the integral satisfy

$$[\frac{\partial h}{\partial x}(x^*(t_f),t_f) - p^*(t_f)]^T \delta x_f + \{g(x^*(t_f),u^*(t_f),a^*,t_f) + \frac{\partial h}{\partial t}(x^*(t_f),t_f) +$$
$$p^*(t_f)[f(x^*(t_f),u^*(t_f),a^*,t_f)]\}^T \delta t_f = 0. \tag{57}$$

Using expanded Hamiltonian function defining as

$$H(x(t),u(t),p(t),a,t) = g(x(t),u(t),a,t) + p^T(t)[f(x(t),u(t),a,t)]. \tag{58}$$

The necessary conditions (7) and (54)-(56) is rewritten by following

$$\dot{x}^*(t) = \frac{\partial H}{\partial p}(x^*(t),u^*(t),p^*(t),a^*,t) \tag{59}$$

$$\dot{p}^*(t) = -\frac{\partial H}{\partial x}(x^*(t),u^*(t),p^*(t),a^*,t) \tag{60}$$

$$0 = \frac{\partial H}{\partial u}(x^*(t),u^*(t),p^*(t),a^*,t) \tag{61}$$

$$0 = \frac{\partial H}{\partial a}(x^*(t),u^*(t),p^*(t),a^*,t) \tag{62}$$

and bounded condition

$$[\frac{\partial h}{\partial x}(x^*(t_f),t_f) - p^*(t_f)]^T \delta x_f + [\frac{\partial h}{\partial t}(x^*(t_f),t_f) + H(x^*(t_f),u^*(t_f),p^*(t_f),a^*,t_f)]\delta t_f = 0. \tag{63}$$

It is clear that if final state is fixed, $\delta t_f$ and $\delta x_f$ are equal to zero, then (63) is vanished and to verify bounded condition $x(t_f) = x_f$ must be satisfied. On the other hand if $t_f$ is fixed and $x(t_f)$ is free, $\delta t_f$ is equal to zero but $\delta x_f$ is arbitrary so (63) is exchanged to

$$\frac{\partial h}{\partial x}(x^*(t_f),t_f) - p^*(t_f) = 0. \tag{64}$$

Suppose $t_f$ is free and $x(t_f)$ is fixed, then $\delta t_f$ is arbitrary but $\delta x_f$ is equal to zero. Therefore (63) is substituted by

$$\frac{\partial h}{\partial t}(x^*(t_f),t_f) + H(x^*(t_f),u^*(t_f),p^*(t_f),a,t_f) = 0. \tag{65}$$

If $t_f$ and $x(t_f)$ are free and independent, to satisfy (63)

$$\frac{\partial h}{\partial x}(x^*(t_f),t_f) - p^*(t_f) = 0 \tag{66}$$

and

$$\frac{\partial h}{\partial t}(x^*(t_f),t_f) + H(x^*(t_f),u^*(t_f),p^*(t_f),a,t_f) = 0 \tag{67}$$

must be satisfied. Now let $t_f$ and $x(t_f)$ are free and related by

$$x(t) = \theta(t) \tag{68}$$

where $\theta(t)$ is a known function, then

$$\delta x_f = [\frac{d\theta}{dt}(t_f)]\delta t_f. \tag{69}$$

By substituting (69) in (63)

$$[\frac{\partial h}{\partial x}(x^*(t_f),t_f) - p^*(t_f)]^T (\frac{d\theta}{dt}(t_f)) + \frac{\partial h}{\partial t}(x^*(t_f),t_f) + H(x^*(t_f),u^*(t_f),p^*(t_f),a^*,t_f) = 0 \tag{70}$$

∎

## Appendix VI

## Proof of Theorem 6

Using transformation equation (15), the cost function is

$$J(u) = H(y(t_0),a,t_0) + \int_{t_0}^{t_f} \{G(y(t),u(t),a,t) + [\frac{\partial H}{\partial y}(y(t),a,t)]^T \dot{y}(t) + \frac{\partial H}{\partial t}(y(t),a,t)\}dt. \tag{71}$$

Note that the first term of $J(u)$ depends on uncertain parameter $a$ and has no efficiency on minimizing $J(u)$, but in order to have min-max solution this term appears in $\delta a$.

Following proves of theorems 5 and 6 conditions 1-8 of theorem 7 are satisfied.